\documentclass[11pt]{article}
\usepackage{amsmath}
\usepackage{amssymb}
\usepackage{amsthm}

\usepackage[margin=2.0cm]{geometry}
\usepackage{hyperref}
\usepackage{bm}

\usepackage{color}
\usepackage{enumerate}
\usepackage{tikz,lmodern}
\usepackage{amsfonts}
\usetikzlibrary{arrows, decorations.markings}

\newtheorem{lem}{Lemma}[section]%
\newtheorem{theorem}[lem]{Theorem}%
\newtheorem{cor}[lem]{Corollary}%
\newtheorem{prob}{Problem}%
\newtheorem{prop}[lem]{Proposition}%

\def\a{\alpha} \def\b{\beta}  \def\d{\delta}

\def\G{\Gamma}

 \def\lg{\langle} \def\rg{\rangle}
\def\nd{\mathrel{\bigm|\kern-.7em/}}

\def\f{\noindent}

\def\PSL{\hbox{\rm PSL}}

\def\PGL{\hbox{\rm PGL}}

\def\Aut{\hbox{\rm Aut\,}}

\def\soc{\hbox{\rm soc}}

\def\Cay{\hbox{\rm Cay}}
\def\Cos{\hbox{\rm Cos }}

\def\mod{\hbox{\rm mod }}

\def\demo{\noindent{\bf Proof}\hskip10pt}

\def\mz{{\mathbb Z}}
\def\qed{\hskip10pt $\Box$\vspace{3mm}}

\renewcommand{\thefootnote}{\fnsymbol{footnote}}

\numberwithin{equation}{section}

\begin{document}
\title{On tetravalent half-arc-transitive graphs}

\author{\\Jin-Xin Zhou$^{\rm a}$\\
{\small\em $^{\rm a}$School of mathematics and statistics, Beijing Jiaotong University, Beijing, 100044, P.R. China}\\}

\renewcommand{\thefootnote}{\fnsymbol{footnote}}
\footnotetext[1]{Corresponding author.
E-mails:
jxzhou@bjtu.edu.cn (J.-X. Zhou)
}

\date{}
\maketitle

\begin{abstract}
Vertex-stabilizers of trivalent edge-transitive graphs have been classified by Tutte, Goldschmidt and some others in several previous papers. Tetravalent half-arc-transitive graphs form an important class of tetravalent edge-transitive graphs. Maru\v si\v c and Nedela (2001) initiated the study of the problem of classifying vertex-stabilizers of tetravalent half-arc-transitive graphs, which has received extensive attention and considerable effort in the literature. 
In this paper, we solve this problem by proving that a group is the vertex-stabilizer of a connected tetravalent half-arc-transitive graph if and only if it is a non-trivial concentric group. Note that a characterization of concentric groups has been given by Maru\v si\v c and Nedela in 2001.
Furthermore, we give an explicit construction of an infinite family of tetravalent half-arc-transitive graphs with automorphism group isomorphic to $A_2^n\wr \mz_2$ and vertex-stabilizers isomorphic to $(D_8^2\times\mz_{2}^{n-6})^2$ for $n\geq7$. These are the first known family of basic tetravalent half-arc-transitive graphs of bi-quasiprimitive type.

\end{abstract}

\section{Introduction}
Throughout this paper, groups are assumed to be finite, and graphs are assumed to be finite, simple, connected and undirected. For the group-theoretic and graph-theoretic terminology not defined here we refer the reader to \cite{Bondy-book,Wielandt-book}.

Let $\G$ be a graph. Then $V(\Gamma)$, $E(\Gamma)$, $A(\G)$ and $\Aut(\Gamma)$  denote its vertex set, edge set, arc set, and full automorphism group, respectively. The vertex stabilizer in $\Aut(\G)$ will be called the {\em vertex stabilizer} of $\G$. We say that $\G$ is {\em vertex-transitive}, {\em edge-transitive} or {\em arc-transitive} if $\Aut(\G)$ is transitive on $V(\G)$, $E(\G)$ or $A(\G)$, respectively. Determining vertex-stabilizers is an active topic in the study of transitive graphs. In 1947, Tutte~\cite{Tutte1947} proved that the vertex-stabilizers of trivalent arc-transitive graphs have order dividing $3\cdot2^4$. In 1980, Goldschmidt~\cite{Goldschmidt1980} obtained an important extension of this result by showing that the vertex-stabilizers of trivalent edge-transitive graphs have order dividing $3\cdot2^7$. For more results on vertex-stabilizers on trivalent edge-transitive graphs, one may see \cite{Djokovic1972,Sims1967,Tutte1959,Weiss}.

For tetravalent edge-transitive graphs, by a result of Gardiner~\cite{Gardiner,GardinerII} (see also \cite[Lemma~2.3]{FangLiXu}), if $\G$ is a connected tetravalent $2$-arc-transitive graph (namely, $\Aut(\G)$ is transitive on the set of the triples $(u, v, w)$ of three pairwise distinct vertices such that $v$ is adjacent to both $u$ and $w$), then the vertex-stabilizers of $\G$ have order dividing $2^4\cdot3^6$. If $\G$ is a connected tetravalent vertex- and edge-transitive but not $2$-arc-transitive graph, then it is easy to prove that the vertex-stabilizer of $\G$ is a $2$-group. In \cite{MN-JGT}, Maru\v si\v c and Nedela initiated the problem of classifying vertex-stabilizers of tetravalent half-arc-transitive graphs.

A graph $\G$ is called {\em half-arc-transitive\/} or {\em HAT} for short if it is vertex- and edge-transitive but not arc-transitive. The study of HAT graphs was initiated by Tutte in \cite{Tutte1966}. Over the last half a century, a fair amount of work have been done on HAT graphs, for which the reader may refer to the survey papers~\cite{CPS2015,M1998}. In the extensive study of HAT graphs, a lot of efforts have been devoted to tetravalent HAT graphs, see for example, \cite{M-DM,MM1999,MN-JGT,MP1999,MX1997,RSparl2019,Spiga-Xia2021,Xia-2021}.

In \cite{MN-JGT}, it was proved that the vertex stabilizer of every tetravalent HAT graph is a concentric group, which is defined as below. A group $H=\lg a_1, \dots, a_n\rg$ is called {\em concentric} if $|\lg a_i, \dots, a_j\rg|=2^{j-i+1}$ for all $1\leq i<j\leq n$ and there exists an isomorphism $\phi:\ \lg a_1, \dots, a_{n-1}\rg\rightarrow \lg a_2, \dots, a_{n}\rg$ such that $a_i^\phi=a_{i+1}$ for $i=1, \dots, n-1.$ The investigation of concentric groups was initiated by Glauberman in \cite{Glauberman1969,Glauberman1971}. Maru\v si\v c and Nedela~\cite{MN-JGT} made an in-depth analysis of such groups, and they proved that these groups are just the point stabilizers of the transitive permutation groups with a non-self-paired suborbit of length $2$. With this result, it was proved in \cite{MN-JGT} that for any given concentric group $H$, there exists a tetravalent graph $\G$ such that $\Aut(\G)$ has a subgroup $G$ which is half-arc-transitive on $\G$ with vertex stabilizer isomorphic to $H$. For the vertex-stabilizers of tetravalent HAT graphs, Spiga commented in \cite{Spiga2016}: ``it is not at all clear whether
every concentric group is the vertex stabilizer of a connected half-arc-transitive graph of valency four." So an intriguing problem is:

\begin{prob}\label{prob}
Determine the vertex-stabilizers of connected tetravalent HAT graphs.
\end{prob}

Despite considerable effort made, the above problem has been wide open so far. One reason perhaps is constructing tetravalent HAT graphs with prescribed vertex stabilizers often turns out to be challenging. In 2005, Maru\v si\v c~\cite{M-DM} proved that every nontrivial elementary abelian $2$-group is the vertex stabilizer of a connected tetravalent HAT graph. Tetravalent HAT graphs with non-abelian vertex-stabilizers are much more elusive. The first know example of such graph was found by Conder and Maru\v si\v c~\cite{C-M-JCTB2003} in 2003. Their graph has vertex stabilizer isomorphic to $D_8$, the dihedral group of order $8$. Twelve years later, Conder, Poto\v cnik and \v Sparl~\cite{CPS2015} constructed two new examples of such graphs, which have vertex stabilizers isomorphic to $D_8$ or $D_8\times \mz_2$, respectively. In 2016, Spiga~\cite{Spiga2016} constructed another two new tetravalent HAT graphs with vertex stabilizers isomorphic to $D_8\times D_8$ or $\mathcal{H}_7$, respectively, where $\mathcal{H}_7=\lg a_1,a_2,\dots,a_7\mid a_i^2=1\ {\rm for}\ i\leq 7, (a_ia_j)^2=1\ {\rm for}\ |i-j|\leq 4, (a_1a_{6})^2=a_3, (a_2a_7)^2=a_4, (a_1a_7)^2=a_{5}\rg$. 
In 2021, Xia in \cite{Xia-2021} made a significant progress in constructing tetravalent HAT graphs with non-abelian vertex stabilizers. He constructed the first infinite family of tetravalent HAT graphs with vertex stabilizers isomorphic to $D_8^2\times\mz_2^m$ for $m\geq1$. In the same year, Spiga and Xia~\cite{Spiga-Xia2021} proved that there exist infinitely many tetravalent HAT graphs with vertex stabilizer isomorphic to $\mathcal{H}_7\times\mz_2$ or $D_8\times\mz_2^m$ for $m\geq0$. They also proved that for any non-trivial concentric group $H$ of order up to $2^8$, there exist infinitely many tetravalent HAT graphs with vertex stabilizer isomorphic to $H$. To the best of our knowledge, the above mentioned groups: $\mz_2^m, D_8\times\mz_2^m, D_8^2\times\mz_2^m$ for $m\geq 1$, $\mathcal{H}_7$ and $\mathcal{H}_7\times\mz_2$, are the only known concentric groups which are vertex-stabilizers of tetravalent HAT graphs.

In this paper, we first prove the following result which shows that for any non-abelian concentric group of order $2^n$ with $n\geq 9$, there exist infinitely many finite connected tetravalent HAT Cayley graphs with vertex-stabilizer isomorphic to $H$. 

\begin{theorem}\label{th:lift}
Let $\G$ be a connected tetravalent graph such $\Aut(\G)$ has a subgroup $G$ which is half-arc-transitive on $\G$. Assume the vertex stabilizer $G_v$ of a vertex $v\in V(\G)$ in $G$ is a non-abelian group of order at least $2^9$.
If $p\geq|G|$ is a prime, then there exists a connected tetravalent half-arc-transitive graph $\Sigma$ such that $\Aut(\Sigma)$ has a normal $p$-subgroup $P$ such that $\Aut(\Sigma)/P\cong G$.
\end{theorem}

This together with some known results enable us to establish the following theorem, which completely solves Problem~\ref{prob}.

\begin{theorem}\label{th:main}
A group is isomorphic to the vertex-stabilizer of a connected tetravalent half-arc-transitive graph if and only if it is a non-trivial concentric group.
\end{theorem}

Let $\G$ be a graph. Assume that $G\leq\Aut(\G)$ is such that $G$ is transitive on the vertex set of $\G$. Let $N$ be a normal subgroup of $G$ such that $N$ is intransitive on $V(\G)$. The {\em normal quotient graph\/} $\G_N$ of $\G$ relative to $N$ is defined as the graph with vertices the orbits of $N$ on $V(\G)$ and with two different orbits adjacent if there exists an edge in $\G$ between the vertices lying in those two orbits. If $\G_N$ and $\G$ have the same valency, then we say that $\G$ is a {\em normal cover} of $\G_N$. An important approach in the study of graphs with prescribed symmetry conditions is to analysis the graphs using normal quotient reduction, and one of remarkable achievements is the characterization of $2$-arc-transitive graphs by Praeger~\cite{Praeger,Praeger1993}. In 2016, this approach was used in \cite{Praeger-normal3} to investigate tetravalent $G$-HAT graphs. Let $\G$ be a connected tetravalent $G$-HAT graph with $G\leq\Aut(\G)$, namely, $G$ is transitive on $V(\G)$ and on $E(\G)$ but intransitive on $A(\G)$. Let $\mathcal{OG}(4)$ denote the family of pairs of $(\G, G)$ for which $\G$ is a connected tetravalent $G$-HAT graph. A pair $(\G, G)$ in $\mathcal{OG}(4)$ is said to be {\em basic} if $\G_N$ has valency at most $2$ for every non-trivial normal subgroup $N$ of $G$. If a pair $(\G, G)$ in $\mathcal{OG}(4)$ is basic, then we say $\G$ is a {\em basic $G$-HAT graph}. A tetravalent basic $\Aut(\G)$-HAT graph $\G$ will be simply called a tetravalent {\em basic HAT graph}.

Following \cite{Praeger-normal3}, tetravalent basic $G$-HAT graphs $\G$ are divided into the following three types:
\begin{itemize}
  \item {\em quasiprimitive type}:\ every non-trivial normal subgroup of $G$ is transitive on $V(\G)$;
  \item {\em bi-quasiprimitive type}:\  every non-trivial normal subgroup of $G$ has at most two orbits on $V(\G)$, and $\Aut(\G)$ has at least one normal subgroup with exactly two orbits on $V(\G)$;
  \item {\em cycle type}: $G$ has at least one non-trivial normal subgroup $N$ such that $\G_N$ is a cycle of length $n\geq3$.
\end{itemize}
Tetravalent basic $G$-HAT graphs of quasiprimitive and bi-quasiprimitive types have been successfully described in \cite{Praeger-normal3,PozPra21}, and those of cycle type have also been extensively studied in \cite{Praeger-normal1,PozPra24,PozPra2x}. For more work on tetravalent basic $G$-HAT graphs, we refer the reader to \cite{Praeger-normal2,PozPra19}. In contrast, however, relatively little is known about tetravalent basic HAT graphs. In \cite{Xia-2021}, Xia constructed an infinite families of tetravalent HAT graphs with automorphism group isomorphic to $A_{2^n}$ for $n\geq7$. To the best knowledge of us, this is the only known infinite family of tetravalent basic HAT graphs which are quasiprimitive. Our next theorem constructs an infinite family of tetravalent basic HAT graphs which are bi-quasiprimitive.

\begin{theorem}\label{th:bi-quasi-type}
There exists a tetravalent basic HAT graph with full automorphism group the wreath product $A_{2^n}\wr\mz_2$ and vertex-stabilizers isomorphic to $(D_8^2\times\mz_{2}^{n-6})^2$ for each $n\geq7$.
\end{theorem}

\section{Preliminaries}
For a positive integer $n$, denote by $\mz_n$ the cyclic group of order $n$, by $D_{2n}$ the dihedral group of order $2n$,
and by $A_n$ and $S_n$ the alternating group and symmetric group of degree $n$, respectively. For a group $G$, we denote by $1$, $Z(G)$ and $\Aut(G)$, the identity element, the center and the automorphism group of $G$, respectively. For two groups $M$ and $N$, $N\rtimes M$ denotes a semidirect product of $N$ by $M$. For a subgroup $H$ of a group $G$, denote by $C_G(H)$ the centralizer of $H$ in $G$ and by $N_G(H)$ the normalizer of $H$ in $G$. Furthermore, the subgroup $\bigcap_{g\in G}H^g$ is called the {\em core} of $H$ in $G$, and if $\bigcap_{g\in G}H^g=1$, then we say that $H$ is {\em core-free} in $G$. If $G$ is a permutation group on a set $\Omega$ and $H$ is a permutation group on a set $\Delta$, then denote $G\wr H$ the wreath product of $G$ by $H$.

Let $G$ be a permutation group on a set $\Omega$. If $G$ fixes a subset $\Delta\subseteq \Omega$ setwise, then denote by $G^\Delta$ the permutation group on $\Delta$ induced by $G$.  Let $\a\in \Omega$. The stabiliser in $G$ of a point $\a\in\Omega$ is
the subgroup $G_\a=\{g\in G \mid \a^g=\a\}$ of $G$. If $G_\a$ is trivial for every $\a\in\Omega$, then we say that $G$ is \emph{semiregular} on $\Omega$, and \emph{regular} if in addition $G$ is transitive on $\Omega$.  If $G$ is transitive on $\Omega$ and the only partitions of $\Omega$ preserved by $G$ are either the singleton subsets or the whole of $\Omega$, then $G$ is said to be {\em primitive} on $\Omega$. If every non-trivial normal subgroup of $G$ is transitive on $\Omega$, then we say that $G$ is {\em quasiprimitive} on $\Omega$. Note that every primitive group is quasiprimitive.

\begin{lem}{\rm\cite[Theorem~2.2]{LMP2009}}\label{lem:primitive2-group}
Let $G$ be a quasiprimitive permutation group a set $\Omega$ of degree $2^n$ with $n$ an integer. Let $N$ be a minimal normal subgroup of $G$. Then $G$ is primitive, and one of the following holds:
\begin{enumerate}[{\rm (1)}]
  \item $G$ is an affine group, $N=\mz_2^n$ and $G\leq{\rm AGL}(n,2)$ with $n\geq2$;
  \item $G$ is almost simple, and $N=T\cong A_{2^n}$ or $\PSL(2, q)$ with $q+1=2^n$ for $n\geq3$, where $q$ is a prime;
  \item $G$ is of production type, and $N=T^\ell$ with $\ell\geq2$, and $T=A_{2^s}$ or $\PSL(2, q)$ with $q+1=2^s$ for $s\geq3$, where $q$ is a prime.
\end{enumerate}
\end{lem}

Let $\G$ be a graph. Then $V(\Gamma)$, $E(\Gamma)$, $D(\G)$ and $\Aut(\Gamma)$  denote its vertex set, edge set, arc set, and full automorphism group, respectively. For $v\in V(\G)$, let $\G(v)$ denote the set of vertices of $\G$ adjacent to $v$. A graph is said to be {\em $k$-valent} if there exists an integer $k$ such that $|\G(v)|=k$ for all vertices $v\in V(\G)$.

Let $G\leq\Aut(\G)$. We say that $\G$ is {\em $G$-vertex-transitive\/}, {\em $G$-edge-transitive\/} or {\em $G$-arc-transitive\/} if $G$ is transitive on $V(\G)$, $E(\G)$ or $A(\G)$, respectively. If $\G$ is $G$-vertex-transitive and $G$-edge-transitive but not $G$-arc-transitive, then we say that $\G$ is {\em $G$-half-arc-transitive}, or {\em $G$-HAT} for short. 
When $G=\Aut(\G)$, a $G$-vertex-transitive, $G$-edge-transitive, {$G$-arc-transitive\/} or {$G$-half-arc-transitive} graph $\G$ is simply called {\em vertex-transitive\/}, {\em edge-transitive\/}, {\em are-transitive\/} or {\em half-arc-transitive\/} (or {\em HAT} for short), respectively.

Let $G$ be a finite group with a core-free subgroup $H$, and an element $g\in G\setminus H$ such that $G=\lg H, g\rg$. The {\em coset graph} $\G:=\Cos(G, H, H\{g, g^{-1}\}H)$ is the graph with vertex set $[G: H]=\{Hx\mid x\in G\}$ such that $Hx$ and $Hy$ are adjacent if and only if $yx^{-1}\in H\{g, g^{-1}\}H$.

The following lemmas are about coset graphs, of which the first can be obtained from \cite{Subdussi}.

\begin{lem}\label{lem:cosetgraph1}
Let $\G$ be a $G$-vertex-transitive and $G$-edge-transitive graph with $G\leq\Aut(\G)$ and let $u\in V(\G)$. Then $\G\cong\Cos(G, G_u, G_u\{g, g^{-1}\}G_u)$, where $g\in G$ and $\{u, u^g\}\in E(\G)$.
\end{lem}

Let $G$ be a finite group with a core-free subgroup $H$ and an element $g\in G\setminus H$. Let $\G=\Cos(G, H, $ $H\{g,$ $ g^{-1}\}H)$.
It is easy to see that $G$ acts faithfully on $V(\G)$ by right multiplication, and this action induces a group of automorphisms of $\G$. For convenience, we shall view $G$ as a subgroup of $\Aut(\G)$.

\begin{lem}{\rm\cite[Lemma~2.1]{Li-Lu-Zhang-JCTB}}\label{lem:cosetgraph2}
Let $G$ be a finite group with a core-free subgroup $H$ and an element $g\in G\setminus H$. Let $\G=\Cos(G, H, H\{g, g^{-1}\}H)$. Then we have
\begin{enumerate}
  \item [{\rm (1)}]\ $\G$ is $G$-vertex-transitive and $G$-edge-transitive;
  \item [{\rm (2)}]\ $\G$ is $G$-arc-transitive if and only if $HgH=Hg^{-1}H$;
  \item [{\rm (3)}]\ $\G$ is connected if and only if $G=\lg H, g\rg$;
  \item [{\rm (4)}]\ the valency of $\G$ is equal either to $|H: H^g\cap H|$ if $HgH=Hg^{-1}H$, or to $2|H: H^g\cap H|$ otherwise.
\end{enumerate}
\end{lem}

Given a finite group $G$ and an inverse-closed subset $S\subseteq G\setminus\{1\}$ (that is, $s^{-1}\in S$ for all $s\in S$), the {\em Cayley graph} $\Cay(G,S)$ of $G$ with respect to $S$ is a graph with vertex set $G$ and edge set $\{\{g,sg\} \mid g\in G,s\in S\}$. Let $\G:=\Cay(G,S)$ be a Cayley graph of a group $G$ with respect to $S$. Then $\G$ is connected if and only if $S$ generates $G$. For any $g\in G$ define
\[
R(g): x\mapsto xg\ \mbox{for $x\in G$ and set $R(G)=\{R(g)\ \mid g\in G\}$.}
\]
Then $R(G)$ is a regular permutation group on $V(\G)$ and is a subgroup of $\Aut(\G)$. Let
\[
\Aut(G,S)=\{\a\in\Aut(G)\mid S^\a=S\}.
\]
It was shown in \cite{Godsil1981} that the normalizer of $R(G)$ in $\Aut(\G)$ is $R(G)\rtimes\Aut(G,S)$. In 1998, Xu introduced the notion of normal Cayley graph. For a group $G$, a Cayley graph $\Cay(G,S)$ of $G$ is called {\em normal} if $R(G)\unlhd \Aut(\Cay(G,S))$ (see \cite{Xu98}).

It is easy to see that if $\Cay(G,S)$ is connected, then $S$ is a generating subset of $G$, and then $\Aut(G,S)$ acts faithfully on $S$. Due to the vertex-transitivity of $\Cay(G, S)$, we have the following lemma.

\begin{lem}\label{lem:cay-faith}
Let $\G=\Cay(G, S)$ be a connected Cayley graph of a group $G$, and let $N=N_{\Aut(\G)}(R(G))$. Then for any $v\in V(\G)$, $N_v$ is faithful on the neighbourhood of $v$ in $\G$.
\end{lem}

\section{Two lemmas}

In this section, we prove two lemmas about $G$-half-arc-transitive graphs.

\begin{lem}\label{lem:iscover}
Let $\G$ be a connected tetravalent $G$-half-arc-transitive graph with $G\leq\Aut(\G)$. Let $N$ be a solvable normal subgroup of $G$. If $G$ is non-solvable, then $\G$ is a normal cover of $\G_N$, and $N$ is semiregular on $V(\G)$.
\end{lem}

\f\demo Consider the normal quotient graph $\G_N$. Take $v\in V(\G)$, and let $O$ be an orbit of $N$ on $V(\G)$ containing $v$. Let $K$ be the kernel of $G$ acting on $V(\G_N)$. Then $K=NK_v$. Since $\G$ is a connected tetravalent $G$-half-arc-transitive graph, $K_v$ is a $2$-group, and since $N$ is solvable, it implies that $K$ is solvable.

If $N$ is transitive on $V(\G)$, then $G=K=NG_v$, and so $G$ is solvable, contrary to the assumption that $G$ is non-solvable. If $N$ has two orbits on $V(\G)$, then $K=NK_v$ is an index $2$ subgroup of $G$. It follows that $G$ is solvable as $K$ is solvable, a contradiction. If $N$ has more than $2$ orbits on $V(\G)$, then $\G_N$ has valency $2$ or $4$. If $\G_N$ has valency $2$, then we may assume that the normal quotient $\G_N\cong C_m$ for some integer $m$. Then $G/K\leq\Aut(\G_N)\cong D_{2m}$, and then $G$ is solvable since $K$ is solvable. Again, we obtain a contradiction. We therefore conclude that $\G_N$ has valency $4$, and $\G$ is a normal cover of $\G_N$. Then $K=N$ is semiregular on $V(\G)$.\hfill\qed

\begin{lem}\footnote{After we proved this lemma, it came to our notice that the same result has been obtained in \cite[Theorem~5]{PP-JCTA} by using a different method.}\label{lem:stab-is-solvable}
Let $\G$ be a connected tetravalent $G$-half-arc-transitive graph with $G\leq\Aut(\G)$. If $G$ is solvable, then the stabilizer $G_v$ is elementary abelian, where $v\in V(\G)$.
\end{lem}

\demo Suppose that $\G$ is a minimum connected tetravalent $G$-half-arc-transitive graph for some $G\leq\Aut(\G)$ such that the vertex stabilizer $G_v$ is non-abelian. We shall seek a contradiction. Let $M$ be a minimal normal subgroup of $G$. Since $G$ is solvable, $M\cong\mz_p^r$ with $p$ a prime and $r\geq 1$. Consider the quotient graph $\G_M$ of $\G$ relative to $M$. Then $\G_M$ has valency $0, 1, 2$ or $4$.

If $\G_M$ has valency $0$, then $M$ is transitive on $V(\G)$. Since $M$ is abelian, $M$ is regular on $V(\G)$, and then $\G$ is a Cayley graph of $M$. In particular, $G_v$ is faithful on $\G(v)$ by Lemma~\ref{lem:cay-faith}. Since $\G$ is $G$-half-arc-transitive, one has $G_v\cong\mz_2$, a contradiction.

If $\G_M$ has valency $1$, then $M$ has two orbits on $V(\G)$. If $M\cong\mz_p^r$ with $p$ an odd prime, then $M$ is semirgular on $V(\G)$, and so $\G$ is a bi-Cayley graph of $M$. Since $M\unlhd G$, by \cite[Lemma~2.3]{HPZ}, $G_v$ is faithful on $\G(v)$, and since $\G$ is $G$-half-arc-transitive, one has $G_v\leq\mz_2\times\mz_2$, a contradiction. If $M$ is a $2$-group, then $G=|V(\G)||G_v|=2|M||G_v|$, and so $G$ is a $2$-group. Since $M$ is a minimal normal subgroup of $G$, it follows that $M$ is contained in the center of $G$, and hence $M\cong\mz_2$. This implies that $\G$ has $4$ vertices, contrary to $\G$ has valency $4$.

If $\G_M$ has valency $4$, then $\G$ is a normal cover of $\G_M$ and $M$ is semiregular on $V(\G)$. Furthermore, $M$ is the kernel of $G$ acting on $V(\G_M)$ and $\G_M$ is a tetravalent $G/M$-half-arc-transitive graph. Let $O$ be the orbit of $G$ on $V(\G)$ containing $v$. Then $G_O=MG_v$. Since $M$ is semiregular on $V(\G)$, the stabilizer of $O\in V(\G_M)$ in $G$ is $G_O/M=G_v/(M\cap G_v)\cong G_v$. So $G_O/M$ is non-abelian. This is impossible by the minimality of $\G$.

Finally assume that $\G_M$ has valency $2$. Then $\G_M\cong C_m$ with $m\geq 3$. Let $K$ be the kernel of $G$ on $V(\G_M)$. Let $B$ be the orbit of $M$ containing $v$. Let $C$ and $D$ be the two orbits of $M$ on $V(\G)$ adjacent to $B$. Since $G$ is edge-transitive on $\G$, every orbit of $M$ does not contain two adjacent vertices. Furthermore, there exists $g\in G$ such that $\{B,C\}^g=\{B,D\}$, and so the subgraphs $\G[C\cup B]$ and $\G[D\cup B]$ of $\G$ induced by $B\cup C$ and $B\cup D$, respectively, are two isomorphic graphs of valency $2$. Let $\G(v)\cap C=\{c_1, c_2\}$ and $\G(v)\cap D=\{d_1, d_2\}$. Since $C$ and $D$ are two orbits of $K$, it follows that $K_v$ setwise fixes $\{c_1, c_2\}$ and $\{d_1, d_2\}$. Assume that $K_v$ pointwise fixes either $\{c_1, c_2\}$ or $\{d_1, d_2\}$. Without loss of generality, assume that $K_v$ pointwise fixes $\{d_1, d_2\}$. Then $K$ is intransitive on the edges between $B$ and $D$. Since $G$ is edge-transitive on $\G$, there exists $g\in G$ such that $\{B, D\}^g=\{B, C\}$. It follows that $K^g=K$ is also intransitive on the edges between $B$ and $C$. This implies that $K_v$ also fixes $c_1$ and $c_2$. Then $K_v$ fixes all neighbours of $v$. By the connectivity of $\G$ as well as the transitivity of $G$ on $V(\G)$, we conclude that $K_v$ will fix all vertices of $\G$, and so $K_v=1$. This implies that $K=M$ is semiregular on $V(\G)$. Since $\G_M\cong C_m$, one has $G/M\leq\Aut(\G_M)\cong D_{2m}$. So $G_v\cong\mz_2$, a contradiction. We therefore conclude that $K_v$ is transitive on  $\{c_1, c_2\}$ and on  $\{d_1, d_2\}$. Since $\G$ is $G$-half-arc-transitive, it follows that $\{c_1,c_2\}$ and $\{d_1,d_2\}$ are also two orbits of $G_v$ on $\G(v)$. Since $\G_M\cong C_m$, one has $G_v=K_v$ and $G/K\cong \mz_m$.

If $M$ is a $2$-group, then $K$ is a $2$-group, and since $G/K$ is cyclic, it follows that $G$ has a normal Sylow $2$-subgroup, say $P$. Since $M$ is minimal normal in $G$, $M$ is contained in the center of $P$. Since $B$ is an orbit of $M$, $K_v$ fixes all vertices of $B$. Let $u\in B$ be adjacent to $d_1$ and $u\neq v$. Since $K_v$ fixes $u$, one has $u$ is adjacent to $d_1$ and $d_2$. So $\{v, u, d_1, d_2\}$ induces a $4$-cycle of $\G$. So $\G$ has girth at most $4$. If $\G$ has girth $3$, then from \cite[Theorem~5.1]{PW-JCTB2007} we conclude that $G_v\leq\mz_2^2$, a contradiction. Now let $\G$ have girth $4$. If $\G$ has two vertices with the same set of neighbours, then by \cite[Lemma~4.3]{PW-JCTB2007}, we have $\G\cong C_{2^\ell}[2K_1]$ and since $\G$ is $G$-half-arc-transitive, one has $G_v\leq\mz_2^\ell$, a contradiction. So we may assume that every two vertices of $\G$ do not have the same set of neighbours. Let $C'=(v, d_1, x, y)$ be a $4$-cycle of $\G$ passing through $\{v, d_1\}$. Then $y\in\{c_1, c_2, d_2\}$. Suppose that $y\in\{c_1, c_2\}$. If $x\in B$, then since $K_v=K_x$, one has $\G(x)=\G(v)$, a contradiction. If $x\notin B$, then $x\notin B\cup C\cup D$, and so $(B, D, x^M, C)$ is a $4$-cycle of $\G_M$. Since $\G_M$ is a cycle, it follows that $\G_M=(B, D, x^M, C)$. Then $|G|=4|K|$ and so $G$ is a $2$-group. Since $M$ is minimal normal in $G$, one has $M\cong\mz_2$ and then $B=\{u,v\}$. So $\G(u)=\G(v)$, a contradiction. Suppose that $y=d_2$. Since $K_{d_1}$ fixes $d_2$, it follows that $\G(d_1)\cap x^M$ and $\G(d_2)\cap x^M$ are both equal to the orbit of $K_{d_1}$ on $x^M$. As every two vertices of $\G$ do not have the same set of neighbours, we have $x^M=B$ and $x=u$. So $C'=(v, d_1, u, d_2)$. This implies that $C'$ is the unique $4$-cycle passing through the edge $\{v, d_1\}$. Since $G$ is edge-transitive on $\G$, in $\G$ there exists only one $4$-cycle passing through every edge. Then in $\G$ there exist exactly two $4$-cycles passing through every vertex.

Now let $\Sigma$ be the graph with vertices the $4$-cycles of $\G$, and two $4$-cycles are adjacent in $\Sigma$ whenever they share a common vertex of $\G$. Then $\Sigma$ has valency $4$, and $G$ acts faithfully on $V(\Sigma)$. Since $\G$ is $G$-half-arc-transitive, $G$ is vertex- and edge-transitive on $\Sigma$. Note that in the $4$-cycle $C'=(v, d_1, u, d_2)$, the $4$ arcs $(v, d_1), (v, d_2), (u, d_1)$ and $(u, d_2)$ are in the same orbit of $G$ on the arc set of $\G$, and the remaining $4$ arcs of $C'$ are in the same orbit of $G$ on the arc set of $\G$. Since $\G$ is $G$-half-arc-transitive, it follows that the setwise stabilizer $G_{C'}$ of $C'$ in $G$ is intransitive on $V(C')$. This implies that $G_{C'}$ is intransitive on $\Sigma({C'})$ and so $\Sigma$ is a tetravalent $G$-half-arc-transitive graph. Since $G_v=K_v=K_u$, one has $G_{C'}=G_v$. This is contrary to the minimality of $\G$.

Next we assume that $M\cong\mz_p^r$ with $p$ an odd prime. Then $K=M\rtimes K_v$. Since $K/M\unlhd G/M$ and $(|K_v|, |M|)=1$, one has
\[G/M=N_{G/M}(K/M)=N_{G/M}(MK_v/M)=N_G(K_v)M/M.\]
This implies that $G=N_G(K_v)M$. Let $g\in G$ be such that $v^g\in\G(v)$. Since $\G$ is $G$-half-arc-transitive, we have $K_v\cap K_{v^g}=K_v\cap K_v^g$ has index $2$ in $K_v$. Since $G=N_G(K_v)M$, one has $g=xy$ for some $x\in N_G(G_v)$ and $y\in M$. So $K_v\cap K_v^g=K_v\cap K_v^y$. Since $|K_v: K_v\cap K_v^y|=2$, one has $y$ has order $p$.

We claim that $K_v\cap K_v^y=C_{K_v}(y)$. Let $h\in K_v\cap K_v^y$. Then $h^{y^{-1}}\in K_v$ and so $h^{-1}h^{y^{-1}}\in K_v$. Since $y\in M$ and $M\unlhd X$, one has $h^{-1}h^{y^{-1}}=h^{-1}yhy^{-1}\in M$. So $h^{-1}h^{y^{-1}}\in K_v\cap M=1$ and so $h=h^{y^{-1}}$. This implies that $h\in C_{K_v}(y)$ and so $K_v\cap K_v^y\leq C_{K_v}(y)$. Clearly, $C_{K_v}(y)\leq K_v\cap K_v^y$, so $K_v\cap K_v^y=C_{K_v}(y)$.

Since $M$ is minimal normal in $G$, there exist $h_2,\dots,h_r\in N_G(K_v)$ such that
\[M=\lg y\rg\times\lg y^{h_2}\rg\times\cdots\times\lg y^{h_r}\rg.\]
Then $C_{K_v}(y^{h_i})=(C_{K_v}(y))^{h_i}<K_v^{h_i}=K_v$ for $i=2,3,\dots,r$. Note that $\bigcap_{i=1}^rC_{K_v}(y^{h_i})\leq C_{K_v}(M)=1$.
Since every $C_{K_v}(y^{h_i})$ is a maximal subgroup of $K_v$, the Frattini subgroup of $K_v$ is trivial. So $K_v (=G_v)$ is elementary abelian. We obtain a contradiction, completing the proof.
\hfill\qed

\section{Vertex-stabilizers of tetravalent half-arc-transitive graphs}

In this section, we shall prove every non-trivial concentric group can be a vertex stabilizer of a tetravalent half-arc-transitive graph. We first prove Theorem~\ref{th:lift} by virtue of a theorem in \cite{Spiga-P-PAMS2019}.


\medskip
\f{\bf Proof of Theorem~\ref{th:lift}}\ Let $\G$ be a connected tetravalent $G$-half-arc-transitive graph with $G\leq\Aut(\G)$. Take $v\in V(\G)$. Assume the stabilizer $G_v$ is a non-abelian group of order at least $2^9$. By Theorem~6 of \cite{Spiga-P-PAMS2019}, there exists a $q$-fold regular cover $\Sigma$ of $\G$ for some power $q$ of $p$, such that the maximal lifted group of automorphisms of $\G$ is $G$. Let $\tilde{G}$ be the lift of $G$. Then $\tilde{G}$ has a normal $p$-subgroup, say $P$
such that $\tilde{G}/P\cong G$. Furthermore, $\Sigma$ is $\tilde{G}$-half-arc-transitive. Let $A=\Aut(\Sigma)$. Then $\tilde{G}\leq A$. To prove this theorem, it suffices to prove that $\tilde{G}=A$

Suppose on the contrary that $\tilde{G}<A$. We shall seek a contradiction. Let $B\leq A$ be such that $\tilde{G}$ is a maximal subgroup of $B$. Then $\tilde{G}=N_B(P)$. Since $G_v$ is non-abelian, by Lemma~\ref{lem:stab-is-solvable}, $G$ is non-solvable, and since $\tilde{G}/P\cong G$, $\tilde{G}$ is non-solvable. It then follows that $B$ is non-solvable.

Let $\tilde{v}$ be a vertex of $\Sigma$. Then $B=\tilde{G}B_{\tilde{v}}$, and then $|B: \tilde{G}|\mid |B_{\tilde{v}}|$. Since $\tilde{G}_{\tilde{v}}\leq B_{\tilde{v}}$ and $|G_v|\geq 2^9$, it follows from a result of Gardiner in \cite{Gardiner,GardinerII}~(see also example \cite[Lemma~2.3]{FangLiXu}) that $B$ is not transitive on the $2$-arcs of $\Sigma$. So $B_{\tilde{v}}$ is a $2$-group. Since $p\geq|G|$ and $|\tilde{G}|=|P||G|$, it follows that $P$ is also a Sylow $p$-subgroup of $B$. By Sylow's theorem, we have $|B: \tilde{G}|=|B: N_B(P)|=kp+1$ for some integer $k\geq 1$.

Let $K$ be the core of $\tilde{G}$ in $B$, and let $\Omega=\{\tilde{G}b\mid b \in B\}$. Then $B$ acts primitively on $\Omega$ by right multiplication, and $K$ is the kernel of this action. So $B/K$ is a primitive permutation group on $\Omega$ with a point stabilizer $(B/K)_\a=\tilde{G}/K$, where $\a=\tilde{G}$. Furthermore, $|\Omega|=|B: \tilde{G}|=2^\ell$ for some $\ell\geq 1$ as $|B: \tilde{G}|\mid |B_{\tilde{v}}|$. Recall that $kp+1=|B: \tilde{G}|=|\Omega|$ with $k\geq 1$. So $|\Omega|>p\geq|G|$.

Let $Q=P\cap K$. Since $P\unlhd G$, $Q$ is the unique Sylow $p$-subgroup of $K$, and since $K\unlhd B$, one has $Q\unlhd B$. Since $\tilde{G}=N_B(P)$, one has $Q<P$. Since $PK/K\cong P/(P\cap K)=P/Q$, it follows that $PK/K$ is a non-trivial normal $p$-subgroup of $\tilde{G}/K$.

By Lemma~\ref{lem:primitive2-group}, $B/K$ is a primitive permutation group on $\Omega$ such that one of parts (1)--(3) of Lemma~\ref{lem:primitive2-group} happens.

If Lemma~\ref{lem:primitive2-group}~(1) happens, then $B/K$ is an affine group. Let $N/K=\soc(B/K)$. Then $N/K\cong\mz_2^\ell$, and $N/K$ is the unique minimal normal subgroup of $B/K$.

Recall that $Q=P\cap K\unlhd B$. Consider the normal quotient graph $\Sigma_Q$. Since $B$ is non-solvable, by Lemma~\ref{lem:iscover}, $\Sigma_Q$ has valency $4$, and $\Sigma$ is a normal cover of $\Sigma_Q$, and $Q$ is semiregular on $V(\Sigma)$


For every subgroup $F$ of $B$ with $Q\leq F$, we let $\bar{F}=F/Q$. Then $\bar{\tilde{G}}/\bar{K}$ is a core-free maximal subgroup of $\bar{B}/\bar{K}$. Furthermore, $\bar{N}/\bar{K}\cong N/K\cong\mz_2^\ell$, and $\bar{N}/\bar{K}$ is the unique minimal normal subgroup of $\bar{B}/\bar{K}$.

Let $\bar{C}=C_{\bar{B}}(\bar{K})$. Since $\bar{P}\cap \bar{K}=1$ and $\bar{P},\bar{K}\unlhd \bar{\tilde{G}}$, it follows that $\bar{P}\leq\bar{C}$, and hence $\bar{C}\bar{K}/\bar{K}$ is a non-trivial normal subgroup of $\bar{B}/\bar{K}$. So $\bar{N}/\bar{K}\unlhd \bar{C}\bar{K}/\bar{K}$. Then $\bar{N}\leq \bar{C}\bar{K}$, and then $\bar{N}=C_{\bar{N}}(\bar{K})\bar{K}$. Let $\bar{E}$ be a Sylow $2$-subgroup of $C_{\bar{N}}(\bar{K})$. As $C_{\bar{N}}(\bar{K})\bar{K}/\bar{K}=\bar{N}/\bar{K}\cong\mz_2^\ell$, it follows that $C_{\bar{N}}(\bar{K})/Z(\bar{K})\cong\mz_2^\ell$. This implies that $\bar{E}$ is normal in $C_{\bar{N}}(\bar{K})$ and $|\bar{E}|\geq 2^\ell$.
Then $\bar{E}$ is characteristic in $C_{\bar{N}}(\bar{K})$, and since $C_{\bar{N}}(\bar{K})\unlhd \bar{N}$, one has $\bar{E}\unlhd \bar{N}$. Now let $\bar{O}$ be the maximal normal $2$-subgroup of $\bar{N}$. Then $1<\bar{E}\leq\bar{O}$. Since $\bar{O}$ is characteristic in $\bar{N}$ and $\bar{N}\unlhd\bar{B}$, one has $\bar{O}\unlhd\bar{B}$. Then $O\unlhd B$ and $|O|=2^r|Q|$ for some $r\geq\ell$. Consider the quotient graph $\Sigma_O$. Since $B$ is non-solvable, by Lemma~\ref{lem:iscover}, $\Sigma$ is a normal cover of $\Sigma_O$ and $O$ is semiregular on $V(\Sigma)$. Hence, $|O|\mid |V(\Sigma)|$. Note that $|V(\Sigma)|=|P||V(\G)|$ and $|O|=2^r|Q|$ with $r\geq \ell$. Since $p\geq |G|$ and $G$ is transitive on $V(\G)$, one has $p\geq|V(\G)|$. This implies that $2^r\mid |V(\G)|$, and then $2^r\mid |G|$. In particular, $2^\ell\leq |G|\leq p$. This is contrary to $2^\ell=|\Omega|>p$.


If Lemma~\ref{lem:primitive2-group}~(2) happens, then $\soc(B/K)\cong A_{2^\ell}$ or $\PSL(2, q)$ with $q+1=2^\ell$, where $\ell\geq3$, $q\equiv 3\ (\mod 4)$ and $q$ is a prime. If $\soc(B/K)\cong A_{2^\ell}$, then the point stabilizer $(B/K)_\a=\tilde{G}/K\cong A_{2^\ell-1}$ or $S_{2^\ell-1}$, and since $\ell\geq3$, the point stabilizer $(B/K)_\a$ does not have a non-trivial normal $p$-subgroup, which is impossible. If $\soc(B/K)\cong \PSL(2, q)$ with $q+1=2^\ell$, then the point stabilizer $(B/K)_\a=\tilde{G}/K\cong \mz_q:\mz_{\frac{q-1}{2}}$ or $\mz_q:\mz_{q-1}$. Since $PK/K$ is a non-trivial normal $p$-subgroup of $\tilde{G}/K$ such that $p\geq |G|=|\tilde{G}/P|$, we have $p=q$ and $PK/K=\mz_q$. Since $\tilde{G}$ is non-solvable, $K$ is non-solvable. It follows that $|K/(K\cap P)|\geq 60$ and hence $|K|\geq 60|K\cap P|$. Note that $|\tilde{G}|=\frac{(q-1)}{k}\cdot |PK/K||K|=\frac{(q-1)}{k}\cdot |P/(P\cap K)||K|\geq \frac{60(q-1)}{k}\cdot |P|$, where $k=1$ or $2$. Since $p\geq |G|=|\tilde{G}: P|$, one has $p\geq 30(p-1)$, a contradiction.

If Lemma~\ref{lem:primitive2-group}~(3) happens, then we may assume that $\Omega=\Delta^k$ for some set $\Delta$ of cardinality $2^{\ell'}$ and integer $k\geq2$ such that $k\cdot 2^{\ell'}=2^\ell$, and that $B/K\leqslant Y=H \wr S_k$ with $\soc(B/K)=\soc(Y)$, and
either $\soc(H)\cong A_{2^{\ell'}}$ or $\soc(H)\cong\PSL(2, q)$ with $q+1=2^{\ell'}$, where $\ell'\geq3$, $q\equiv 3\ (\mod 4)$ and $q$ is a prime. Now let $T=\soc(H)$ and take $\d\in \Delta$. Then $H_\d$ is maximal in $H$ and $T_\d\neq 1$. Furthermore, $\soc(Y)=\soc(B/K)=T^k$. Let $S=\soc(B/K)$ and take $\omega=(\d,\d,\ldots,\d)\in\Delta^k=\Omega$. Then $S_\omega=(T_\d)^k$. As $\tilde{G}/K$ is the stabilizer of the point $\a=\tilde{G}$ of $\Omega$ in $B/K$, we have $\tilde{G}/K=((B/K)_\omega)^g$ for some $g\in B/K$. Without loss of generality, we may assume that $\tilde{G}/K=(B/K)_\omega$.

Assume that $S=\{(g_1, g_2, \ldots, g_k)\mid g_i\in T\}$. Let
\[\begin{array}{l}
T_1=\{(g,1,\ldots,1) \mid g\in T\}, T_2=\{(1,g,1,\ldots,1) \mid g\in T\}, \ldots, T_k=\{(1,1,\ldots,1,g) \mid g\in T\}.
\end{array}
\]
Then $S=T_1\times T_2\times\cdots\times T_k$, and $S_\omega=(T_1)_\omega\times (T_2)_\omega\times\cdots\times (T_k)_\omega$. Since $B/K$ is primitive on $\Omega$, $(B/K)_\omega$ acts transitively on $\{T_1, T_2,\dots, T_k\}$ by conjugation, and the kernel of this action is contained in $H^k$.
Similarly, let
\[\begin{array}{l}
H_1=\{(g,1,\ldots,1) \mid g\in H\}, H_2=\{(1,g,1,\ldots,1) \mid g\in H\}, \ldots, H_k=\{(1,1,\ldots,1,g) \mid g\in H\}.
\end{array}
\]
Then $H^k=H_1\times H_2\times\cdots \times H_k$. Clearly, every $(T_i)_\omega\neq 1$, and $(B/K)_\omega$ also acts transitively on $\{(T_1)_\omega, (T_2)_\omega, \dots, (T_k)_\omega\}$ by conjugation. Let $L$ be the kernel of $(B/K)_\omega$ acting on $\{(T_1)_\omega, (T_2)_\omega, \dots, (T_k)_\omega\}$. Then $L$ fixes every $(T_i)_\omega$ and hence fixes every $T_i$. This implies that $L\leq (H^k)\cap (B/K)_\omega$, and so $L\leq (H_1)_\omega\times (H_2)_\omega\times\dots\times (H_k)_\omega$. Recall that either $\soc(H)\cong A_{2^{\ell'}}$ or $\soc(H)\cong \PSL(2, q)$ with $q+1=2^{\ell'}$. If $\soc(H)\cong A_{2^{\ell'}}$, then $A_{2^{\ell'}-1}\leq (H_i)_\omega\leq S_{2^{\ell'}-1}$ and $(T_i)_\omega\cong A_{2^{\ell'}-1}$ for each $1\leq i\leq k$. If $\soc(H)\cong \PSL(2, q)$, then $\mz_q:\mz_{\frac{q-1}{2}}\leq (H_i)_\omega\leq \mz_q:\mz_{q-1}$ and $(T_i)_\omega\cong \mz_q:\mz_{\frac{q-1}{2}}$ for each $1\leq i\leq k$. This implies that the centralizer of $S_\omega$ in $(H_1)_\omega\times (H_2)_\omega\times\dots\times (H_k)_\omega$ is trivial. So $C_L(S_\omega)=1$ and hence $C_{(B/K)_\omega}(S_\omega)=1$.

Since $S=\soc(B/K)\unlhd B/K$, one has $S_\omega\unlhd (B/K)_\omega=\tilde{G}/K$. Since $1\neq PK/K\unlhd \tilde{G}/K$, if $(PK/K)\cap S_\omega=1$, then $PK/K\leq C_{(B/K)_\omega}(S_\omega)=1$, a contradiction. Thus, $1\neq PK/K\cap S_\omega \unlhd S_\omega$. This implies that for each $1\leq i\leq k$, $T_i\cong \PSL(2, q)$ and $(T_i)_\omega\cong\mz_q\rtimes \mz_{\frac{q-1}{2}}$, where $q+1=2^{\ell'}$, $\ell'\geq3$, and $q\equiv 3\ (\mod 4)$ and $q$ is a prime. Furthermore, $(H_i)_\omega\cong\mz_q\rtimes \mz_{\frac{q-1}{2}}$ or $\mz_q\rtimes \mz_{{q-1}}$ for $1\leq i\leq k$. Since $1\neq PK/K\cap S_\omega \unlhd S_\omega$, one has $p\mid |S_\omega|$.
Note that $|\tilde{G}/K|=|G||P|/|K|$ and
\[(T_1)_\omega\times (T_2)_\omega\times\cdots\times (T_k)_\omega\unlhd (B/K)_\omega=\tilde{G}/K.\]
As $p\geq |\tilde{G}/P|=|G|$, one has $p$ is the largest prime divisor of $|\tilde{G}|$. Since $p\mid |S_\omega|$, one has $p=q$. Since $p\geq |\tilde{G}/P|=|G|$, one has $p\geq (\frac{q-1}{2})^k\geq (\frac{q-1}{2})^2=\frac{(p-1)^2}{4}$. So $4p\geq p^2-2p+1$, forcing $p\leq 5$. However, since $G$ is non-solvable, one has $|G|\geq60$, and then $p\geq|G|\geq60$, a contradiction. This completes the proof.\hfill\qed

In the following of this section, we shall use the following notation.
Let $n\geq 2$ be an integer, and let $H$ be a concentric group of order $2^n$. By definition (see also \cite[Theorem~5.5]{MN-JGT}), $H$ is generated by $n$ involutions, say $a_1, a_2, \dots, a_n$, such that $|\lg a_i,\dots, a_j\rg|=2^{j-i+1}$ for all $1\leq i<j\leq n$. Furthermore, there exits a rotary group isomorphism
\begin{equation}\label{eq-phi}
\phi: \lg a_1, a_2, \dots, a_{n-1}\rg\rightarrow\lg a_2, a_3, \dots, a_n\rg
\end{equation}
such that $a_i^\phi=a_{i+1}$ for $i=1, \dots, n-1$.
Let $B=\lg a_1, a_2, \dots, a_{n-1}\rg$ and $C=\lg a_2, a_3, \dots, a_n\rg.$

Let $h\in B$, and define a permutation $\tau_h$ on $H$ as follows:
\begin{equation}\label{eq-th}
\tau_h: b^{\tau_h}=b^\phi, (a_mb)^{\tau_h}=a_1hb^\phi\ {\rm for}\ b \in B.
\end{equation}

\begin{lem}\label{lem:R(H)concentric}
For any $h\in B$, $\tau_h^{-1}R(a_i)\tau_h=R(a_{i+1})$ for $i=1,2,\dots,n-1$.
\end{lem}

\f\demo Let $i=1, 2, \dots, n-1$ and $c\in C$. Then $b^\phi=c$ for some $b\in B$. Then
\[c^{\tau_h^{-1}R(a_i)\tau_h}=(b^{\phi})^{\tau_h^{-1}R(a_i)\tau_h}=(b^{\phi\tau_h^{-1}})^{R(a_i)\tau_h}=b^{R(a_i)\tau_h}=(ba_i)^{\tau_h}=b^\phi a_i^\phi=b^\phi a_{i+1}=c^{R(a_{i+1})},\]
and
\[\begin{array}{ll}
(a_1hc)^{\tau_h^{-1}R(a_i)\tau_h}&=(a_1hb^{\phi})^{\tau_h^{-1}R(a_i)\tau_h}=(a_mb)^{R(a_i)\tau_h}=(a_mba_i)^{\tau_h}\\
&=a_1h(ba_i)^{\phi}=a_1hb^\phi a_i^\phi=a_1hca_{i+1}=(a_1hc)^{R(a_{i+1})}.
\end{array}\]\hfill\qed

\begin{lem}\label{lem:2-creator}
Let $G=\lg\tau_h, R(H)\rg$. If $H$ is non-abelian, then $R(H)$ is core-free in $G$, and the coset graph $\G=\Cos(G, R(H), R(H)\{\tau_h^{-1}, \tau_h\}R(H))$ is a connected tetravalent $G$-half-arc-transitive Cayley graph of $G_1$, and the stabilizer of the vertex $R(H)$ of $\G$ in $G$ is $R(H)$, where $G_1$ is the stabilizer of $1\in H$ in $G$.
\end{lem}

\f\demo It is easy to see that $G$ is a transitive permutation group on $H$ and $\tau_h\in G_1$. Let $a=\tau_h^{-1}$ and $b=R(a_1)a$.
Then $R(a_1)=ab^{-1}$. By Lemma~\ref{lem:R(H)concentric}, $R(a_j)=a^{j}b^{-1}a^{-j+1}=a^{j-1}R(a_1)a^{-(j-1)}$ for $j=1,2,\dots,n$.
It then follows that $G=\lg R(a_1), a\rg=\lg R(a_1), b\rg=\lg a,b\rg$. Since $R: g\mapsto R(g)$ for all $g\in H$ is an isomorphism from $H$ to $R(H)$, it follows that $|\lg R(a_i),\dots, R(a_j)\rg|=2^{j-i+1}$ for all $1\leq i\leq j\leq n$. Since $H$ is non-abelian, by \cite[Corollary~5.3]{MN-JGT}, $H$ is not normal in $G$ and then $(G, H)$ is a $\vec{2}$-creator (see \cite[p.23]{MN-JGT} for the definition of the $\vec{2}$-creator). By \cite[Theorem~3.1]{MN-JGT}, $R(H)$ is core-free in $G$, and $R(H)\tau_hR(H)\neq R(H)\tau_h^{-1}R(H)$ since $R(H)$ is non-abelian. Note that $R(H)^{\tau_h}\cap R(H)=\lg R(a_2),\dots, R(a_n)\rg$. It follows that the coset graph $\Cos(G, R(H), R(H)\{\tau_h^{-1}, \tau_h\}R(H))$ is a connected tetravalent $G$-half-arc-transitive graph. Since $R(H)$ is core-free, the stabilizer of the vertex $R(H)$ of $\G$ in $G$ is $R(H)\cong H$. Since $R(H)$ is regular on $H$, one has $R(H)\cap G_1=1$. This implies that $\G$ is a Cayley graph of $G_1$.\hfill\qed

\begin{cor}\label{cor-Cayley}
Let $n$ be a positive integer and $H$ a non-abelian concentric group of order $2^n$ with $n\geq 9$. Then there exist infinitely many finite connected tetravalent half-arc-transitive Cayley graphs with vertex-stabilizer isomorphic to $H$.
\end{cor}

\f\demo Combining Lemma~\ref{th:lift} with Lemma~\ref{lem:2-creator}, we can obtain this corollary.\hfill\qed

\begin{prop}\label{cor:concentric-stabilizer}
Let $n$ be a positive integer and $H$ a concentric group of order $2^n$. Then there exist infinitely many finite connected tetravalent HAT graphs with vertex-stabilizer isomorphic to $H$.
\end{prop}

\f\demo If $H$ is abelian, then by \cite[Theorem~5.5]{MN-JGT}, $H$ is an elementary abelian $2$-group, and combining \cite[Theorem~1.1]{M-DM} and \cite[Theorem~1.1]{Spiga-Xia2021}, there exist infinitely many finite connected tetravalent half-arc-transitive graphs with vertex-stabilizer isomorphic to $H$.

Assume that $H$ is non-abelian. If $|H|\leq 2^8$, then by \cite[Theorem~1.4]{Spiga-Xia2021}, there exist infinitely many finite connected tetravalent half-arc-transitive graphs with vertex-stabilizer isomorphic to $H$. If $|H|\geq 2^9$, then by Lemma~\ref{lem:2-creator} or \cite[Corollary~7.4]{MN-JGT}, there exists a finite connected tetravalent graph $\G$ and a subgroup $G\leq\Aut(\G)$ such that $G$ is half-arc-transitive on $\G$ with vertex-stabilizer $G_v\cong H$ for some vertex $v$ of $\G$. By Theorem~\ref{th:lift}, there exist infinitely many finite connected tetravalent half-arc-transitive graph with vertex-stabilizer isomorphic to $H$.\hfill\qed

\f{\bf Proof of Theorem~\ref{th:main}}\ By \cite[Theorem~1.1]{MN-JGT}, we may see that the vertex-stabilizer of a  connected tetravalent half-arc-transitive graph is a concentric group. By Proposition~\ref{cor:concentric-stabilizer}, every non-trivial concentric group is isomorphic to the vertex-stabilizer of some connected tetravalent half-arc-transitive graph.\hfill\qed

\section{Tetravalent basic $G$-half-arc-transitive graphs with non-abelian vertex stabilizers}

\subsection{A reduction}

The following lemma is easy to prove.

\begin{lem}\label{cor:cover-basic}
Let $\G$ be a connected tetravalent $G$-half-arc-transitive graph. Then $G$ has a normal subgroup $N$ such that $\G_N$ a tetravalent basic $G/N$-half-arc-transitive graph.
\end{lem}

\f\demo 
Let $N\unlhd G$ be maximal subject to the condition that $\G$ is a normal cover of the normal quotient graph $\G_N$. Then $\G_N$ is a connected tetravalent $G/N$-half-arc-transitive graph, and $N$ is semiregular on $V(\G)$. Furthermore, $\G_N$ is $G/N$-basic. \hfill\qed

The following lemma is about tetravalent basic $G$-half-arc-transitive graph with non-abelian vertex-stabilizers.

\begin{lem}\label{th:basic}
Let $\G$ be a connected tetravalent basic $G$-half-arc-transitive graph with non-abelian vertex-stabilizers. Then the socle of $G$ is a non-abelian minimal normal subgroup of $G$.
\end{lem}

\f\demo\ Take $v\in V(\G)$ and let $H=G_v$. Then $H$ is a non-abelian concentric group of order $2^n$ for some integer $n\geq 3$. By Lemma~\ref{lem:stab-is-solvable}, $G$ is non-solvable. Let $M$ be a minimal normal subgroup of $G$. Since $\G$ is tetravalent basic $G$-half-arc-transitive graph, the normal quotient $\G_M$ has valency $0, 1$ or $2$.

Suppose that $M$ is abelian. If $M$ is transitive on $V(\G)$, then $G=MH$ and then $G$ is solvable, a contradiction.
If $M$ has two orbits, say $B_0, B_1$, on $V(\G)$, then $\{B_0, B_1\}$ is a $G$-invariant partition of $V(\G)$. Since $G$ is transitive on $V(\G)$, the setwise stabilizer $G_{B_0}$ of $B_0$ in $G$ is $MH$, which is an index $2$ subgroup of $G$. Since $MH$ is solvable, it follows that $G$ is solvable, a contradiction. If $M$ has more than $2$ orbits on $V(\G)$, then $\G_M\cong C_m$ for some integer $m\geq3$. Let $K$ be the kernel of $G$ acting on $V(\G)$. Then $K=MK_v$ and $G/K\leq \Aut(\G_M)\cong D_{2m}$. Noticing that $K=MK_v$ is solvable, it follows that $G$ is solvable, a contradiction.

We therefore conclude that $M$ is non-abelian. To prove this theorem, it suffices to prove that $M$ is just the socle $\soc(G)$ of $G$. Suppose on the contrary that $R$ is another minimal normal subgroup of $G$. The arguments in the above paragraph implies that $R$ is also non-abelian. Furthermore, $M\times R\unlhd G$, and $R$ is contained in the centralizer of $M$ in $G$. If $M$ is transitive on $V(\G)$, then $G=MH$, and then $G/M$ is a $2$-group. Since $R\cong RM/M\leq G/M$, $R$ would be a $2$-group. This is impossible because $R$ is a non-abelian minimal normal subgroup of $G$. If $M$ has two orbits, say $B_0, B_1$, on $V(\G)$, then $\{B_0, B_1\}$ is a $G$-invariant partition of $V(\G)$. Since $G$ is transitive on $V(\G)$, the setwise stabilizer $G_{B_0}$ of $B_0$ in $G$ is $MH$, which is an index $2$ subgroup of $G$. Since $R$ is a non-abelian minimal normal subgroup of $G$, $R$ has no index two subgroups. It follows that $R\leq MH$, and hence $R\cong RM/M\leq MH/M$. This implies that $R$ would be a $2$-group. Again, this is impossible because $R$ is a non-abelian minimal normal subgroup of $G$. If $M$ has more than $2$ orbits on $V(\G)$, then $\G_M\cong C_m$ for some integer $m\geq3$. Let $K$ be the kernel of $G$ acting on $V(\G_M)$. Then $K=MK_v$ and $G/K\leq \Aut(\G_M)\cong D_{2m}$. Recall that $R$ is a minimal normal subgroup of $G$. If $R\cap K\neq 1$, then $R\leq K$, and then $R\cong RM/M\leq K/M\leq K_v$. This forces that $R$ is solvable, a contradiction. If $R\cap K=1$, then $R\leq G/K\leq D_{2m}$, and then $R$ would be solvable, a contradiction.\hfill\qed

\subsection{Constructing tetravalent basic $G$-half-arc-transitive graphs}

Next we shall give a construction of tetravalent basic $G$-half-arc-transitive graphs. Let $W$ be a simple primitive permutation group on the set $\Delta=\{1,2,\dots,2^n\}$ with $n\geq3$. Suppose that $W$ has an element $a$ fixing $1$ and a regular subgroup $H$ satisfying the following conditions:
\begin{enumerate}
  \item [{\rm (C1)}]\ $W=\lg H,a\rg$,
  \item [{\rm (C2)}]\ $H$ is generated by $n$ involutions: $h_1, h_2, \dots, h_n$,
  \item [{\rm (C3)}]\ $|\lg h_i, \dots, h_j\rg|=2^{j-i+1}$ for all $1\leq i<j\leq n$,
  \item [{\rm (C4)}]\ $h_i^a=h_{i+1}$ for $i=1,2,\dots,n-1$.
\end{enumerate}

Let $m\geq 1$ be an integer. For $0\leq i\leq m-1$, let $\Delta_i=\{k+i\cdot 2^n\mid k\in\Delta\}$. Let $\Omega=\bigcup_{i=0}^{m-1}\Delta_i$.
Let $\tau$ be a permutation on $\Omega$ such that
\[\tau: i\mapsto 2^n+i\mapsto 2\cdot 2^n+i\mapsto 3\cdot 2^n+i\mapsto \cdots\mapsto (m-1)\cdot2^n+i\mapsto i, \ {\rm for}\ i\in \Delta.\]
Let $H_0=H$, and let $H_i=H^{\tau^i}$ for $i=1,2,\dots,m-1$.

\begin{theorem}\label{lem:generated}
Let $G=\lg W,\tau\rg$ and $K=H_0\times H_1\times \cdots\times H_{m-1}$. Then $\soc(G)=W\times W^{\tau}\times \cdots \times W^{\tau^{m-1}}$ is a minimal normal subgroup of $G$, and the coset graph $\Cos(G, K, K\{a\tau,(a\tau)^{-1}\}K)$ is a connected tetravalent $G$-half-arc-transitive graph with vertex-stabilizer isomorphic to $K$.
\end{theorem}

\f\demo Note that $G=(W\times W^{\tau}\times \cdots \times W^{\tau^{m-1}})\rtimes\lg\tau\rg$. This implies that $\soc(G)=W\times W^{\tau}\times \cdots \times W^{\tau^{m-1}}$ is a minimal normal subgroup of $G$. Since $K$ is solvable, $K$ is core-free in $G$. Let $\G=\Cos(G, K, K\{a\tau,(a\tau)^{-1}\}K)$. Then $V(\G)=\{Kg\mid g\in G\}$, and $G$ acts faithfully on $V(\G)$ by right multiplication. It follows that $K$ is the vertex-stabilizer of the vertex $K\in V(\G)$ in $G$.

Take $\a\in \Delta_0$. Suppose that $\a^a=\b$ for some $\b\in\Delta_0$. Then $\a^{(a\tau)^m}=\b$. This implies that the restriction of $(a\tau)^m$ on $\Delta_0$ is equal to $a$. So $\lg (a\tau)^m, H_0\rg^{\Delta_0}=W$. Since $1\neq H_0\leq W$ and $W$ is non-abelian simple, the normal closure of $H_0$ in $\lg (a\tau)^m, H_0\rg$ is equal to $W$. So $W\leq \lg K, a\tau\rg$, and then $a\in \lg K, a\tau\rg$. It follows that $\tau\in \lg K, a\tau\rg$ and $W\times W^{\tau}\times \cdots \times W^{\tau^{m-1}}\leq \lg K, a\tau\rg$. So $G=\lg K, a\tau\rg$. This shows that $\G$ is connected.

Note that $H_i=H_0^{\tau^i}$ for $i=1,2,\dots,m-1$. It is easy to verify that
\[\begin{array}{ll}
a\tau:& h_1^{\tau^{m-1}}\mapsto h_1\mapsto h_2^\tau\mapsto  h_2^{\tau^2}\mapsto \dots\mapsto h_2^{\tau^{m-1}}\mapsto h_2\mapsto h_3^\tau\mapsto  h_3^{\tau^2}\mapsto \dots\mapsto h_3^{\tau^{m-1}}\mapsto\\
&\cdots\mapsto h_{n-1}\mapsto h_n^\tau\mapsto  h_n^{\tau^2}\mapsto \dots\mapsto h_n^{\tau^{m-1}}\mapsto h_n.
\end{array}\]
Let $B=\lg h_1, \dots, h_{n-1}\rg\times H_1\times\cdots\times H_{m-1}$, and $C=H_0\times H_1\times\cdots\times H_{m-2}\times \lg h_2^{\tau^{m-1}}, \dots, h_{n}^{\tau^{m-1}}\rg$. Then $B^{a\tau}=C$ and then $K^{a\tau}\cap K=C$. So $|K: K^{a\tau}\cap K|=2$. If $Ka\tau K=K(a\tau)^{-1}K$, then $\G$ will be a cycle. This is clearly impossible because $G$ is non-solvable. Thus, $Ka\tau K\neq K(a\tau)^{-1}K$, and by Lemma~\ref{lem:cosetgraph2}, $\G$ is a connected tetravalent $G$-half-arc-transitive graph. This completes the proof. \hfill\qed

\f{\bf Remark~1.}\ Let $W=A_{2^n}$. If $n\geq7$, then by \cite{Xia-2021}, $W$ has an element $a\in W_1$ and a regular subgroup $H\cong D_8^2\times\mz_2^{n-6}$ satisfying the above conditions (C1)--(C4). If $n\geq 4$, then by \cite{Spiga-Xia2021}, $W$ has an element $a$ and a regular subgroup $H\cong D_8\times\mz_2^{n-3}$ satisfying the above conditions (C1)--(C4). Based on this, we can use Theorem~\ref{lem:generated} to give an explicit construction of a tetravalent basic $G$-half-arc-transitive graph with $\soc(G)=A_{2^n}^m$ for all $m\geq 1$.

\subsection{A family of tetravalent basic half-arc-transitive graphs}
In this section, we shall use the main result in \cite{Xia-2021} to construct a new family of tetravalent basic half-arc-transitive graphs. To do this, we need the following result.

\begin{lem}\label{lem:xia}{\rm\cite[Theorem~1.2]{Xia-2021}}
Let $W=A_{2^n}$ with $n\geq7$. Then $W$ has an element $a\in W_1$ and a regular subgroup $H\cong D_8^2\times\mz_2^{n-6}$ satisfying the above conditions (C1)--(C4). Furthermore, letting $S=\{a, hah_1, a^{-1}, (hah_1)^{-1}\}$ with $1^{ha}=1^{h_1}$, we have $\lg S\rg=W_1$ and $\Aut(W_1, S)=1$.
\end{lem}

The following is the main result of this subsection.

\begin{theorem}\label{th:new-hat-graph}
In Theorem~\ref{lem:generated}, if $m=2$ and $W, a$ and $H$ are the same as in Lemma~\ref{lem:xia}, then the coset graph $\G=\Cos(G, K, K\{a\tau,(a\tau)^{-1}\}K)$ is a connected tetravalent basic half-arc-transitive graph with automorphism group $G$. Furthermore, $\G$ is bi-quasiprimitive.
\end{theorem}

\f\demo In Theorem~\ref{lem:generated}, we let $m=2$. Then $G=(W\times W^\tau)\rtimes\lg\tau\rg$. Let $\G=\Cos(G, K, K\{a\tau,(a\tau)^{-1}\}K)$. We claim that $G=\Aut(\G)$. Suppose that $G<\Aut(\G)$. Then we may choose a subgroup $M\leq\Aut(\G)$ such that $G$ is maximal in $M$.

Let $N$ be the core of $G$ in $M$, and let $\Delta$ be the set of right cosets of $G$ in $M$. Then $M$ primitively acts on $\Delta$ by right multiplication, and $N$ is the kernel of this action. Let $\d=G$. Then the stabilizer of $\d$ in $M$ is $G$.

Suppose first that $N=1$. Then $M$ acts faithfully on $\Delta$. So we may view $M$ as a primitive permutation group on $\Delta$. Note that $|\Delta|=|M: G|\mid |M: V(\G)|$. It follows that $|\Delta|=2^m$ for some integer $m\geq 3$.
Now apply Lemma~\ref{lem:primitive2-group} to $M$ and $\soc(M)$.

If Lemma~\ref{lem:primitive2-group}~(1) happens, then $\soc(M)\cong\mz_2^m$, and $M=\soc(M)\rtimes G$ and $G$ is an irreducible subgroup of ${\rm GL}(m,2)$. Note that $G=(A_{2^n}\times A_{2^n})\rtimes\mz_2$. In particular, $A_{2^n}\times A_{2^n}$ is a subgroup of ${\rm GL}(m,2)$. By \cite[Proposition~5.5.7]{Kleidman-Liebeck}, we have $m\geq 2m_1$, where $m_1$ is the minimal dimension of a non-trivial irreducible linear representation of $A_{2^n}$ in characteristic $2$. By \cite[Proposition~5.3.7]{Kleidman-Liebeck}, we have $m_1=2^n-2$. It follows that $m\geq 2^{n+1}-4$, and then $2^m\geq 2^{2^{n+1}-4}$. On the other hand, by Lemma~\ref{lem:iscover}, $2^m=|\soc(M)|\mid |V(\G)|=|G: K|=\frac{2\cdot(2^n!)^2}{4\cdot (2^n)^2}$. Note that the largest power of $2$ dividing $2^n!$ is $2^n-1$. So $2^m\mid \frac{2\cdot(2^{2^n-1})^2}{4\cdot (2^n)^2}=2^{2^{n+1}-2n-3}$. Since $2^m\geq 2^{2^{n+1}-4}$, one has $2^{n+1}-4\leq2^{n+1}-2n-3$ and then $2n\leq1$, a contradiction.

If Lemma~\ref{lem:primitive2-group}~(2) happens, then $\soc(M)\cong A_{2^m}$ or $\PSL(2,q)$. For the former, $M\cong A_{2^m}$ or $S_{2^m}$, and then $M_\d\cong A_{2^m-1}$ or $S_{2^m-1}$. For the latter, $M\cong\PSL(2,q)$ or $\PGL(2,q)$, and $M_\d\cong \mz_q:\mz_{\frac{q-1}{2}}$ or $\mz_q:\mz_{{q-1}}$. Both cases are impossible because $M_\d=G=(A_{2^n}\times A_{2^n})\rtimes\mz_2$.

If Lemma~\ref{lem:primitive2-group}~(3) happens, then $\soc(M)_\d=A_{2^m-1}^2$ or $(\mz_q:\mz_{\frac{q-1}{2}})^2$. Since $\soc(M)_\d\unlhd M_\d=G$ and $\soc(G)=A_{2^n}\times A_{2^n}$, we have $A_{2^n}\times A_{2^n}\unlhd \soc(M)_\d$, which is impossible. 

We now suppose that $N\neq1$. Then $A_{2^n}\times A_{2^n}=\soc(G)\leq N$. If $N=\soc(G)$, then $M/N$ is a primitive permutation group on $\Delta$ with $G/N\cong\mz_2$ as the stabilizer of $\d$ in $M/N$. So $M/N\cong D_{2p}$ for some odd prime $p$ and $|\Delta|=p$, contrary to $|\Delta|=2^m$. So $N>\soc(G)$ and then $N=G\unlhd M$. Then $|M: G|=2$.

Let $C=C_M(\soc(G))$. Then $C\cap \soc(G)=1$ and $M/C\leq\Aut(\soc(G))=(S_{2^n}\times S_{2^n})\rtimes\mz_2$. Since $|M: \soc(G)|=4$ and $C_G(\soc(G))=1$, one has $|C|=1$ or $2$.

For $i=0,1$, since $H_i$ is regular on $\Delta_i$, we may identify $\Delta_i$ with $H_i$, and identify $1$ and $1+2^n$, respectively, with the identities of $H_0$ and $H_1$. Then $H_i$ acts on $\Delta_i$ by rightly multiplication, and $G$ has a primitive product action on $\Delta_0\times \Delta_1$. Let $X$ be the stabilizer of $(1, 1+2^n)$ in $G$. Then $X=(W_1\times W_1^\tau)\rtimes\lg\tau\rg$, which is maximal in $G$. As $K=H_0\times H_1$ is regular on $\Delta_0\times \Delta_1$, one has $X\cap K=1$ and $G=XK$. This implies that $X$ is regular on $V(\G)$ and then $\G$ is a Cayley grapn of $X$. Since $h_1^\tau\notin K\cap K^{a\tau}=H_0\times\lg h_2^\tau, \dots, h_n^\tau\rg$, one has $Ka\tau h_1^\tau\neq Ka\tau$. So $\{Ka\tau, Ka\tau h_1^\tau\}$ is an orbit of $K$ acting on the neighbourhood of $K$ in $\G$. Since $\G$ is a tetravalent $G$-half-arc-transitive graph, it follows that $\{K(a\tau)^{-1}, K(a\tau h_1^\tau)^{-1}\}$ is another orbit of $K$ acting on the neighbourhood of $K$ in $\G$.

Note that $a\tau$ fixes $(1, 1+2^n)$. So $a\tau\in X$. Let $h\in H$ be such that $1^{ha}=1^{h_1}$. Then $ha\tau h_1^\tau$ also fixes $(1, 1+2^n)$. So $ha\tau h_1^\tau\in X$. Since $Kha\tau h_1^\tau=Ka\tau h_1^\tau$, one has $\G=\Cay(X, T)$ with
\[T=\{a\tau, ha\tau h_1^\tau, (a\tau)^{-1}, (ha\tau h_1^\tau)^{-1}\}.\]
Since $X$ is regular on $V(\G)$, we may identify $X$ with $V(\G)$. Then $M$ is also a permutation group on $X$, and $X$ acts regularly on $X$ by right multiplication.

We first claim that $2\mid |N_M(X)\cap \Aut(X, T)|$. This is clearly true if $|C|=2$. Assume that $|C|=1$. Then $M\leq\Aut(\soc(G))=(S_{2^n}\times S_{2^n})\rtimes\mz_2$. Let $x=(2, 3)\in {\rm Sym}(\Delta_0)$. Then $W\rtimes\lg x\rg={\rm Sym}(\Delta_0)$, and then $\Aut(\soc(G))=((W\rtimes\lg x\rg)\times (W^\tau\rtimes\lg x^\tau\rg))\rtimes\lg\tau\rg$. So $\Aut(\soc(G))=\soc(G)\rtimes (\lg x, x^\tau\rg\rtimes\lg\tau\rg)$. Then \[M=M\cap \Aut(\soc(G))=(\soc(G)\rtimes\lg \tau\rg)(M\cap \lg x,x^\tau\rg).\]
Since $|M: \soc(G)|=4$, one has $M\cap \lg x,x^\tau\rg=\lg xx^\tau\rg$ and then $M=(W\times W^\tau)\rtimes\lg xx^\tau,\tau\rg$. Since $X=(W_1\times W_1^\tau)\rtimes\lg\tau\rg$, one has $X\rtimes\lg xx^\tau\rg\leq N_M(X)$, as claimed.

Now we know that $N_M(X)\cap \Aut(X, T)$ has an involution, say $g$. Since $g$ normalizes $X$, we have $g$ acts on $X$ by conjugation.
We may view $g$ as an inner automorphism of $M$ induced by $g$. If $g=\tau$, then
$T^\tau=T$, and then
\[\{\tau a, \tau hah_1, a^{-1}\tau, (hah_1)^{-1}\tau\}=\{a\tau, ha\tau h_1^\tau, (a\tau)^{-1}, (ha\tau h_1^\tau)^{-1}\}.\]
If $a\tau=\tau a$, then $\tau a\tau=a\in W\cap W^\tau=1$, a contradiction. If $a\tau=\tau hah_1$, then $\tau a\tau=hah_1\in W\cap W^\tau=1$, a contradiction. If $a\tau=a^{-1}\tau$ or $(hah_1)^{-1}\tau$, then we would have $HaH=Ha^{-1}H$. By Lemma~\ref{lem:cosetgraph2}, $\Cos(W, H, HaH)$ would be a cycle, which is impossible by Theorem~\ref{lem:generated}. We therefore conclude that $g\neq\tau$.

Let $T_2=\{(a\tau)^2, (ha\tau h_1^\tau)^2, (a\tau)^{-2}, (ha\tau h_1^\tau)^{-2}\}$. Then
\begin{equation*}\label{eq-T2}
T_2=\{aa^\tau, (hah_1)(hah_1)^\tau, (aa^\tau)^{-1}, ((hah_1)(hah_1)^\tau)^{-1}\}.
\end{equation*}
Since $g$ acts on $X$ by conjugation, it follows that $g$ setwise fixes $T_2$. Since $\tau$ is an involution, $\tau$ commutes with every element in $T_2$.

Recall that either $|C|=2$ and $M=G\times C$, or $C=1$ and $M=(W\times W^\tau)\rtimes\lg xx^\tau,\tau\rg$ with $x=(2,3)$. So $g=z\tau^i$, where $i\in\{0, 1\}$, $1\neq z\in M$ and $z$ normalizes both $W$ and $W^\tau$. Since $\tau$ commutes with every element in $T_2$, it follows that $z$ setwise fixes $T_2$. Since $z$ normalizes both $W$ and $W^\tau$, it follows that $z$ setwise fixes both $S=\{a, hah_1, a^{-1}, (hah_1)^{-1}\}$ and $S^\tau=\{a^\tau, (hah_1)^\tau, (a^{-1})^\tau, ((hah_1)^{-1})^\tau\}$ by conjugation. However, by Lemma~\ref{lem:xia}, we have $\Aut(W, S)=\Aut(W^\tau, S^\tau)=1$, a contradiction.

By now we have shown that $\Aut(\G)=G$. Note that $\soc(G)=W\times W^\tau$ is the unique minimal normal subgroup. Since $K\leq\soc(G)$, $\soc(G)$ has two orbits on $V(\G)$. This implies that $\G$ is bi-quasiprimitive.
\hfill\qed

\section*{Acknowledgments}
This work was supported by the National Natural Science Foundation of China (12425111, 12071023, 12331013, 12161141005).
	
	

\end{document}